\documentclass[10pt]{amsart}

\usepackage{amsmath, amsthm, amscd, amsfonts, amssymb, graphicx, color}
\usepackage{latexsym}
\usepackage{amscd}
\usepackage{amsthm}
\usepackage{amssymb}
\usepackage{enumerate}
\usepackage{mathabx}

\theoremstyle{plain}
\newtheorem{theorem}{Theorem}
\newtheorem{corollary}[theorem]{Corollary}

\theoremstyle{remark}
\theoremstyle{definition}

\begin{document}

\title{Aichinger equation on commutative semigroups}
\author{J.~M.~Almira}

\subjclass[2010]{}


\address{Departamento de Ingeniería y Tecnología de Computadores, Universidad de Murcia, Murcia, Spain}
\email{jmalmira@um.es}

\maketitle
\begin{abstract}

We consider Aichinger's equation $$f(x_1+\cdots+x_{m+1})=\sum_{i=1}^{m+1}g_i(x_1,x_2,\cdots, \widehat{x_i},\cdots, x_{m+1})$$ for functions defined on commutative semigroups which take values on commutative groups. The solutions of this equation are, under very mild hypotheses, generalized polynomials. We use the canonical form of generalized polynomials to prove that compositions and products of generalized polynomials are again generalized polynomials and that the bounds for the degrees are, in this new context, the natural ones. In some cases, we also show that a polynomial function defined on a semigroup can uniquely be extended to a polynomial function defined on a larger group. For example, if $f$ solves Aichinger's equation under the additional restriction that $x_1,\cdots,x_{m+1}\in \mathbb{R}_+^p$, then there exists a unique polynomial function $F$ defined on $\mathbb{R}^p$ such that $F_{|\mathbb{R}_+^p}=f$. In particular, if $f$ is also bounded on a set $A\subseteq \mathbb{R}_+^p$ with positive Lebesgue measure then its unique polynomial extension $F$ is an ordinary polynomial of $p$ variables with total degree $\leq m$, and the functions $g_i$ are also restrictions to $\mathbb{R}_+^{pm}$ of ordinary polynomials of total degree $\leq m$ defined on $\mathbb{R}^{pm}$. 
\end{abstract}

\markboth{J.~M.~Almira}{Aichinger equation on commutative semigroups}

\section{Introduction}

Compositions and products of generalized polynomials defined on abelian groups are again generalized polynomials, and the proper bounds for their degrees are known. This is not trivial (although it may be not impossible) if you deal with a definition based on Fréchet's unmixed differences functional equation. Things become easier if we consider Fréchet's mixed differences functional equation, which is, in quite general cases, equivalent to the unmixed equation (see, e.g.,  \cite{A01}, \cite{A02}, \cite{A03}, \cite{A04}, \cite{A05},  \cite{A06}, \cite{Djokovich}, \cite{Laczkovich},  \cite{L_Dj}). The proof recently appeared in a paper by Aichinger and Moosbauer \cite{AM}.  In their paper, the authors  introduced the following wonderful functional equation, that characterizes generalized polynomials of degree $\leq m$ on abelian groups:
\begin{equation}\label{A}
f(x_1+\cdots+x_{m+1})=\sum_{i=1}^{m+1}g_i(x_1,x_2,\cdots, \widehat{x_i},\cdots, x_{m+1}),  
\end{equation}
where $\widehat{x_i}$ means that the function $g_i$ does not depend on $x_i$, and they used this equation to prove that the composition $f\circ g$ of the generalized polynomials (defined on abelian groups) $f$, $g$ is a generalized polynomial and $\deg(f\circ g)\leq \deg(f)\cdot \deg(g)$. This same result, with a different proof that does not use \eqref{A}, had already been proved by Leibman in \cite{L}, but the authors of \cite{AM} were unaware of that paper.  Moreover in his paper Leibman also proved that if we consider a composition of several polynomial functions $f_i: G_{i-1}\to G_i$, $i=1,2,...,k$, and $G_k$ is nilpotent (with no other extra hypotheses on the first groups $G_0,G_1,...,G_{k-1}$, that may be non-commutative and quite general, indeed), the composition $f_k\circ f_{k-1}\circ \cdots \circ f_1$ is also a polynomial function, but no precise estimation of its degree was given. Indeed, the question if compositions of polynomial functions defined on arbitrary groups are again polynomial functions was posed by Leibman in \cite{L} and is still open.

 The new equation \eqref{A} was -as the authors confirmed me- proposed by Aichinger and motivated by some algebraic problems that are solved in \cite{M}. I was delighted since I find this equation so nice and natural (e.g., note that $(x+y+z)^2=[2yz+z^2]+[2xz+x^2]+[2xy+y^2]=g_1(y,z)+g_2(x,z)+g_3(x,y)$ and, in general, if you take $f(x)=x^m$ evaluated on $x_1+....+x_{m+1}$, you get a sum of monomials, each one depending on at most m variables) but I never had seen it anywhere before. On the other hand, it is also clear that under strong regularity assumptions  on  $f$, ordinary polynomials are the only solutions of Aichinger's equation. For example, if we assume that $f:\mathbb{R}\to\mathbb{R}$ belongs to $\mathbf{C}^{(m+1)}(\mathbb{R})$, $g_i:\mathbb{R}^m\to\mathbb{R}$ belongs to $\mathbf{C}^{(m+1)}(\mathbb{R}^m)$ for $i=1,\cdots,m+1$,  and these functions satisfy \eqref{A}, then $G(x_1,\cdots,x_{m+1})=f(x_1+\cdots+x_{m+1})$ satisfies that $$\frac{\partial^{m+1}G}{\partial x_1\partial x_2\cdots \partial x_{m+1}}=f^{(m+1)}(x_1+\cdots+x_{m+1})=0$$
which implies that $f^{(m)}=0$, and henceforth $f(t)=a_0+a_1t+\cdots+a_mt^m$ is an ordinary polynomial of degree $\leq m$ .  

A similar result can be demonstrated for smooth functions defined on $\mathbb{R}^p$ for each $p>1$. In that case, the space of solutions is  $$\Pi_{m,\text{tot}}^p=\{\sum_{i_1+\cdots+i_p\leq m}a_{i_1,\cdots,i_p}t_1^{i_1}t_2^{i_2}\cdots t_{p}^{i_p}:  a_{i_1,\cdots,i_p}\in\mathbb{R} \text{ for all } (i_1,\cdots, i_p)\in\mathbb{N}^p\},$$
the set of ordinary polynomials in $p$ variables with total degree $\leq m$. We prove this result as a corollary of the following surprising theorem (for the proof, see  \cite[Theorem 14]{PS}): 
\begin{theorem}[Prager and Schwaiger] \label{PrS} If  $K$ is a field and $ f : K^p\to K $ is an ordinary polynomial separately in each variable (which means that for any $1\leq k\leq p$ and any point $(a_1,\cdots,a_{k-1},a_{k+1},\cdots,a_p)\in K^{p-1}$, the function $f(a_1,\cdots,a_{k-1},x_k,a_{k+1},\cdots,a_p)$ is an ordinary algebraic polynomial in $x_k$) then $f$ is an ordinary polynomial in $p$ variables provided that $K$ is finite or uncountable. Furthermore, for every countable infinite field $K$ there exists a function $f:K^2\to K$ which is an ordinary polynomial separately in each variable and is not a generalized polynomial in both variables jointly. 
\end{theorem} 

In fact, if we assume that $f:\mathbb{R}^p\to\mathbb{R}$ belongs to $\mathbf{C}^{(m+1)}(\mathbb{R}^p)$ and satisfies \eqref{A} with $x_j=(x_{1j},\cdots,x_{pj})\in   \mathbb{R}^p$, $j=1,\cdots,m+1$, for certain smooth functions $g_i\in \mathbf{C}^{(m+1)}(\mathbb{R}^{pm})$, $i=1,\cdots,m+1$,  then  for any $1\leq k\leq p$ and any point $(a_1,\cdots,a_{k-1},a_{k+1},\cdots,a_p)\in \mathbb{R}^{p-1}$, the function $h(t)=f(a_1,\cdots,a_{k-1},t,a_{k+1},\cdots,a_p)$ is an ordinary algebraic polynomial in $t$ of degree at most $m$ since 

\begin{eqnarray*}
&\ & h(t_1+\cdots+t_{m+1}) = f(a_1,\cdots,a_{k-1},t_1+\cdots+t_{m+1},a_{k+1},\cdots,a_p) \\
&=& f\left(\sum_{i=1}^{m+1}(\frac{a_1}{m+1},\cdots,\frac{a_{k-1}}{m+1},t_i,\frac{a_{k+1}}{m+1},\cdots,\frac{a_p}{m+1})\right)\\
&=&\sum_{i=1}^{m+1}\frak{g}_i(t_1,t_2,\cdots, \widehat{t_i},\cdots, t_{m+1}),
\end{eqnarray*}
where
\begin{eqnarray*}
&\ & \frak{g}_i(t_1,t_2,\cdots, \widehat{t_i},\cdots, t_{m+1}) \\
&\ & =g_i(\frac{a_1}{m+1},\cdots,\frac{a_{k-1}}{m+1},t_1,\frac{a_{k+1}}{m+1},\cdots,\frac{a_p}{m+1}, \cdots,
\frac{a_1}{m+1},\cdots,\frac{a_{k-1}}{m+1},t_{i-1},\frac{a_{k+1}}{m+1},\cdots,\frac{a_p}{m+1},\\
&\ & \ \  \frac{a_1}{m+1},\cdots,\frac{a_{k-1}}{m+1},t_{i+1},\frac{a_{k+1}}{m+1},\cdots,\frac{a_p}{m+1}, \cdots, \frac{a_1}{m+1},\cdots,\frac{a_{k-1}}{m+1},t_{m+1},\frac{a_{k+1}}{m+1},\cdots,\frac{a_p}{m+1}) 
\end{eqnarray*} 
which means that $h$ and $\frak{g}_i$, $i=1,\cdots,m+1$  (that are smooth functions) solve Aichinger's equation in the one-dimensional setting.  Then, using  Theorem \ref{PrS}, we conclude that $f$ is an ordinary polynomial of $p$ real variables. Once this is known, it is easy to check that the total degree of $f$ is at most $m$. 

It follows that Aichinger's equation is of interest when the regularity assumptions on $f$ are impossible (because we do not have any topology) or are weak. 

In this paper we consider Aichinger's equation for functions defined on commutative semigroups which take values on commutative groups. In section 2 we introduce the canonical form of generalized polynomials and we use it to prove that a polynomial function defined on a cancellative semigroup can uniquely be extended to a polynomial function defined on a larger group. For example, if $f$ solves Aichinger's equation under the additional restriction that $x_1,\cdots,x_{m+1}\in \mathbb{R}_+^p$, then there exists a unique polynomial function $F$ defined on $\mathbb{R}^p$ such that $F_{|\mathbb{R}_+^p}=f$. Here, $\mathbb{R}_+=]0,\infty [$ is the set of strictly positive real numbers. Finally, we devote section 3 to prove that compositions and products of generalized polynomials are again generalized polynomials and that the bounds for the degrees are, in this new context, the natural ones.

Just to fix some notation, we recall that $(S,+)$ is a semigroup if $+:S\times S\to S$ given by $+(a,b)=a+b$ is an associative binary law, which means that $a+(b+c)=(a+b)+c$ for all $a,b,c\in S$. Moreover, $S$ is unital if there exists $0\in S$ such that $0+a=a+0=a$ for all $a\in S$, it is abelian (also named commutative) if $a+b=b+a$ for all $a,b\in S$, and it is cancellative if $a+c=b+c$ implies $a=b$ for all $a,b,c\in S$. It is well known that every cancellative semigroup $S$ can be extended to a group $G$ such that $G=S-S=\{x-y:x,y\in S\}$ and that, if $S$ is commutative, $G=S-S$ is also commutative. We say that $G=S-S$ is a natural extension of the semigroup $S$ (see \cite[p. 34]{CP}).      

\section{Aichinger equation on commutative semigroups}
Given $(S,+)$ a commutative semigroup and $(H,+)$ a commutative group, we say that $f:S\to H$ is a polynomial function (also named generalized polynomial) of degree $\leq m-1$ if $f$ solves Fréchet's mixed differences functional equation:
$$
\Delta_{h_1}\Delta_{h_2}\cdots \Delta_{h_m}f(x)=0\text{ for all } h_1,\cdots, h_m,x\in G.
$$

A simple statement and proof of the characterization of generalized polynomials as solutions of Aichinger's equation is as follows: 

\begin{theorem} \label{main} Let $(S,+)$ be an abelian semigroup, and $(H,+)$ be an abelian group. Let $f:S\to H$ be a map and let $m\in\{2,3,\cdots\}$. Consider the claims:
\begin{itemize}
\item[$(a)$]  There are functions $g_i:S^{m-1}\to H$, $i=1,2,\cdots, m$ such that $f$ satisfies  
\begin{equation} \label{Aeq}
f(x_1+x_2+\cdots+x_{m})=\sum_{i=1}^{m}g_i(x_1,x_2,\cdots,\widehat{x_{i}},\cdots,x_{m})
\end{equation}
for all  $x_1,\cdots,x_m\in S$, where $\widehat{x_{i}}$ means that $g_i$  does not depend on $x_i$. 
\item[$(b)$] $f$ satisfies:
\begin{equation} \label{FrechSem} \Delta_{h_1}\Delta_{h_2}\cdots \Delta_{h_{m}}f(x_1+x_2+\cdots+x_m)=0
\end{equation} 
for all  $x_1,\cdots,x_m,h_1,\cdots,h_m\in S$. 
\item[$(c)$] $f$ satisfies Fréchet's mixed functional equation:
\begin{equation} \label{Frech}  \Delta_{h_1}\Delta_{h_2}\cdots \Delta_{h_{m}}f(x)=0
\end{equation}  
for all  $x,h_1,\cdots,h_m\in S$.
\end{itemize}
Then $(a)\Rightarrow (b)$. Furthermore:
\begin{itemize}
\item If $S+S=S$, then $(a)\Rightarrow (c)$.
\item If $S+S=S$ and $0\in S$, then $(a)\Leftrightarrow (b)\Leftrightarrow (c)$.
\end{itemize}
\end{theorem} 
\noindent \textbf{Proof. } We prove $(a)\Rightarrow (b)$ by induction on $m$. For $m=1$, \eqref{Aeq} means that $f(x_1)$ is a constant and, of course, this implies that $\Delta_{h_1}f(x_1)=0$ for all $x_1,h_1\in S$. Assume that $(a)\Rightarrow (b)$ holds true for $m-1$ and that $f$ satisties \eqref{Aeq}, which means that 
\[
f(x_1+x_2+\cdots+x_{m})=\sum_{i=1}^{m}g_{i0}(x_1,x_2,\cdots,\widehat{x_{i}},\cdots,x_{m})
\]
for certain functions $g_{i0}$, $i=1,\cdots,m$. 
Let $h_1\in S$ and define $F_1(x)=f(x+h_1)-f(x)$. Then 
\begin{eqnarray*}
F_1(x_1+x_2+\cdots+x_{m})&=& \Delta_{h_1}f(x_1+\cdots+x_m)\\
&=&  f((x_1+h_1)+x_2+\cdots+x_m)-f(x_1+x_2+\cdots+x_m)\\
&=& g_{10}(x_2,x_3,\cdots,x_m)+\sum_{i=2}^{m}g_{i0}(x_1+h_1,x_2,\cdots,\widehat{x_{i}},\cdots,x_{m}) \\
&\ & - g_{10}(x_2,x_3,\cdots,x_m) -\sum_{i=2}^{m}g_{i0}(x_1,x_2,\cdots,\widehat{x_{i}},\cdots,x_{m})\\
&=& \sum_{i=2}^{m}(g_{i0}(x_1+h_1,x_2,\cdots,\widehat{x_{i}},\cdots,x_{m})-g_{i0}(x_1,x_2,\cdots,\widehat{x_{i}},\cdots,x_{m}))\\
&=& \sum_{i=2}^{m}g_{i1}(x_1,x_2,\cdots,\widehat{x_{i}},\cdots,x_{m})\\
\end{eqnarray*}
where 
\[
g_{i1}(x_1,x_2,\cdots,\widehat{x_{i}},\cdots,x_{m})=g_{i0}(x_1+h_1,x_2,\cdots,\widehat{x_{i}},\cdots,x_{m})-g_{i0}(x_1,x_2,\cdots,\widehat{x_{i}},\cdots,x_{m}), i=2,3,\cdots,m.
\]
We can repeat the argument as follows: take $h_2\in S$ and define $F_2(x)=\Delta_{h_2}F_1(x)=F_(x+h_2)-F_1(x)$. Then 
\begin{eqnarray*}
F_2(x_1+x_2+\cdots+x_{m})&=& \Delta_{h_2}F_1(x_1+\cdots+x_m) \\
&=& \Delta_{h_2}\Delta_{h_1}f(x_1+\cdots+x_m) \\
&=&  F_1(x_1+(x_2+h_2)+\cdots+x_m)-F_1(x_1+x_2+\cdots+x_m)\\
&=& g_{21}(x_1,x_3,\cdots,x_m)+\sum_{i=3}^{m}g_{i1}(x_1,x_2+h_2,\cdots,\widehat{x_{i}},\cdots,x_{m}) \\
&\ & - g_{21}(x_1,x_3,\cdots,x_m) -\sum_{i=3}^{m}g_{i1}(x_1,x_2,\cdots,\widehat{x_{i}},\cdots,x_{m})\\
&=& \sum_{i=3}^{m}(g_{i1}(x_1,x_2+h_2,\cdots,\widehat{x_{i}},\cdots,x_{m})-g_{i1}(x_1,x_2,\cdots,\widehat{x_{i}},\cdots,x_{m}))\\
&=& \sum_{i=3}^{m}g_{i2}(x_1,x_2,\cdots,\widehat{x_{i}},\cdots,x_{m})\\
\end{eqnarray*}
where 
\[
g_{i2}(x_1,x_2,\cdots,\widehat{x_{i}},\cdots,x_{m})=g_{i1}(x_1,x_2+h_2,\cdots,\widehat{x_{i}},\cdots,x_{m})-g_{i1}(x_1,x_2,\cdots,\widehat{x_{i}},\cdots,x_{m}), i=3,\cdots,m.
\]
If we repeat this argument $m$ times, we obtain that 
\[
\Delta_{h_m}\Delta_{h_{m-1}}\cdots \Delta_{h_{1}}f(x_1+x_2+\cdots+x_m)=0
 \]
for all  $x_1,\cdots,x_m,h_1,\cdots,h_m\in S$, and $(b)$ holds true. This proves $(a)\Rightarrow (b)$. If $S+S=S$, then $(b)\Rightarrow (c)$ is trivial and, henceforth, $(a)\Rightarrow (c)$ also holds true. 

Let us now assume $(c)$. Then for each $z\in S$ we have that 
$$0 = \Delta_{x_1}\Delta_{x_2}\cdots\Delta_{ x_{m}}f(z) = \sum_{\varepsilon_1,\varepsilon_2,\cdots,\varepsilon_{m}=0}^1(-1)^{m-(\varepsilon_1+\varepsilon_2+\cdots+\varepsilon_m)}f(z+\varepsilon_1x_1+\varepsilon_2 x_2+\cdots +\varepsilon_{m}x_{m}).$$
Now, if we consider $z$ as a constant, there is only one term, $f(z+x_1+x_2+\cdots+x_{m})$, in the sum that appears at the third member of this formula, that depends on all variables $x_1,\cdots,x_{m}$. Thus, adding the reciprocal of $f(z+x_1+x_2+\cdots+x_{m})$ to both sides of the equality and taking reciprocals again, we get an expression of the form 
$$f(z+x_1+\cdots+x_{m})= \sum_{i=1}^{m}g_i(z,x_1,x_2,\cdots,\widehat{x_{i}},\cdots,x_{m})$$
This means that, if $f$ is a polynomial function of degree $\leq m-1$ defined on the commutative semigroup $S$, then all the translations $(\tau_zf)(x)=f(x+z)$ with $z\in S$ satisfy Aichinger's equation \eqref{Aeq}. In particular, when $z=0\in S$ we get that $f$ satisfies Aichinger's equation. This proves $(c)\Rightarrow (a)$ when $0\in S$. $\hfill{ \Box}$

\section{Polynomial functions and an extension theorem}

Given $(S,+)$ a commutative semigroup and $(H,+)$ a commutative group, a map $A:S^k\to H$ is $k$-additive if 
\begin{eqnarray*}
A^k(x_1,\cdots,x_{s-1},x+y,x_{s+1},\cdots,x_k) &=&A^k(x_1,\cdots,x_{s-1},x,x_{s+1},\cdots,x_k) \\
&\ & \ + A^k(x_1,\cdots,x_{s-1},y,x_{s+1},\cdots,x_k)
\end{eqnarray*}
for all $s\in\{1,\cdots,k\}$ and $x_1,\cdots,x_k,x,y\in S$. In other words, $A^k$ is $k$-additive if it is additive in each one of its variables. These maps are natural generalizations of additive functions and they are named multiadditive functions as soon as $k>1$. We abuse of notation and call $0$-additive map to any constant function $A_0:S\to H$. Finally, $A:S^k\to H$ is symmetric if $$A(x_1,\cdots,x_k)=A(x_{\sigma(1)},\cdots, x_{\sigma(k)})$$ for every permutation $\sigma$ of 
$\{1,\cdots,k\}$. For example, $A(x,y)=x+y$  and $B(x,y,z)=xy+xz+yz$, where $x,y,z\in \mathbb{R}$, are symmetric maps. 

Associated to any multiadditive function $A:S^k\to H$, we can consider its diagonalization $diag(A)(x)=A(x,\cdots,x)$, which is a map defined on $S$. It is well known that if multiplication by $m!$ is biyective on $H$ then every polynomial function $f:S\to H$ has a unique representation of the form 
\begin{equation} \label{sumamonomios}
f(x)=A_0+diag(A^1)(x)+\cdots +diag(A^m)(x)
\end{equation} 
where $A^k:S^k\to H$ is a symmetric $k$-additive function for $k=1,\cdots, m$ and $m$ is the degree of the polynomial function $f$.  Moreover, every function of the form \eqref{sumamonomios} is a polynomial function of degree at most $m$. The term $A_k(x)=diag(A^k)(x)$ is called monomial (of degree $k$) . Hence, under these hypotheses,  every polynomial of degree $\leq m$ is a sum of monomials of degrees $\leq m$. 
We include a draft of the proof of these results, just for the sake of completeness. The main idea of the proof is to use the following technical result (whose proof can be found, e.g., in \cite{czerwik} for the case that $S,H$ are $\mathbb{Q}$-vector spaces, but that can be easily adapted to the case that $S$ is a commutative semigroup and $H$ is a commutative group):
\begin{theorem}[Polarization formula] \label{polar} Let $A:S^m\to H$ be a symmetric $m$-additive function, and let $A^*(x)=diag(A)(x)$ be its diagonalization. Then for all $x,h_1,\cdots,h_k\in S$ we have that:
\begin{equation}\label{polarization}
\Delta_{h_1}\Delta_{h_2} \cdots \Delta_{h_k}A^*(x)=\left\{ \begin{array}{llll} 0 & \text{ if } k>m\\ m! A(h_1,\cdots,h_m) & \text{ if } k=m  \end{array}\right . 
\end{equation}
In particular, if $H$ allows division by $m!$, we have that 
\begin{equation}\label{monomialEq}
\frac{1}{m!}\Delta_h^mA^*(x)=A^*(h) \text{ for all } x,h\in S.
\end{equation}
\end{theorem}

Note that  \eqref{polarization} implies that, if multiplication by $m!$ is either surjective in $S$ or inyective in $H$, any symmetric $m$-additive function $A$ is completely determined by its diagonalization $diag(A)$ since $m! A(h_1,\cdots,h_m) =  A(m!h_1,\cdots,h_m)$. In particular, $diag(A)=0$ implies $A=0$. 

The following theorem is well known: 
\begin{theorem}[Canonical representation of polynomial functions] \label{repre} Let $S$ be a commutative semigroup and $H$ be a commutative group. Assume that multiplication by $m!$ is biyective on $H$. Let $f:S\to H$ be a polynomial function of degree at most $m$. Then there exist $k$-additive symmetric functions $A^k: S^k\to H$, $k=0,1,2,\cdots, m$, such that  
 \begin{equation} \label{rep}
 f(x)=A^0+\text{diag}(A^1)(x)+\cdots+\text{diag}(A^m)(x)\ \text{ for all } x\in S
 \end{equation} 
 Moreover, all these functions are polynomial functions of degree $\leq m$ and the representation \eqref{rep} is unique.
\end{theorem}

\noindent \textbf{Proof.} The second claim follows from the polarization formula. To prove the first claim, we proceed by induction on $m$. For $m=0$ is clear that $\Delta_{x_1}f=0$ implies that $f(x)=A^0$ for a certain constant $A^0\in H$. Assume the result holds true for polynomial functions of degree $\leq m-1$ and let $f:S\to H$ be a polynomial function of degree at most $m$. Then $\Delta_{x_1}\Delta_{x_2}\cdots \Delta_{x_m}f(x)$ does not depend on $x$ since 
$\Delta_{x_{m+1}}\Delta_{x_1}\Delta_{x_2}\cdots \Delta_{x_m}f(x)=0$ for all $x_{m+1},x\in S$. Hence 
$$A(x_1,\cdots,x_m)=\Delta_{x_1}\Delta_{x_2}\cdots \Delta_{x_m}f(x)$$
is well defined as a function $A:S^m\to H$, and is symmetric since the operators $\Delta_{x_k}$ are pairwise 
commuting. Moreover, $A$ is $m$-additive since $\Delta_{x+y}=\Delta_x\Delta_y+\Delta_x+\Delta_y$ implies that
\begin{eqnarray*}
A(x_1+y_1,x_2,\cdots,x_m) &=& \Delta_{x_1+y_1}\Delta_{x_2}\cdots \Delta_{x_m}f\\
&=& \Delta_{y_1}\Delta_{x_1}\Delta_{x_2}\cdots \Delta_{x_m}f+ \Delta_{x_1}\Delta_{x_2}\cdots \Delta_{x_m}f+  \Delta_{y_1}\Delta_{x_2}\cdots \Delta_{x_m}f\\
&=&  \Delta_{x_1}\Delta_{x_2}\cdots \Delta_{x_m}f+  \Delta_{y_1}\Delta_{x_2}\cdots \Delta_{x_m}f\\
&=& A(x_1,x_2,\cdots,x_m)+A(y_1,x_2,\cdots,x_m).
\end{eqnarray*}
and the symmetry gives the additivity with respect to the other variables. 
Our assumption that multiplication by $m!$ is biyective on $H$ implies that the operation $y\to \frac{1}{m!}y$ is well defined on $H$. Thus, we can define $A^m(x_1,\cdots,x_m)=\frac{1}{m!}A(x_1,\cdots,x_m)$ -which is also $m$-additive- and $f_m(x)=diag(A^m)(x)$. The polarization formula implies that 
\[
\Delta_{x_1}\Delta_{x_2}\cdots \Delta_{x_m}f_m(x)=m!A^m(x_1,\cdots,x_m)=A(x_1,\cdots,x_m)=\Delta_{x_1}\Delta_{x_2}\cdots \Delta_{x_m}f(x)
\]
Hence 
\[
\Delta_{x_1}\Delta_{x_2}\cdots \Delta_{x_m}(f-f_m)(x)=0
\]
and we can apply the induction hypothesis (note that, if the operation $y\to \frac{1}{m!}y$ is well defined on $H$, the same holds with the operation  $y\to \frac{1}{k!}y= (\frac{m!}{k!})\frac{1}{m!}y$ for $1\leq k\leq m-1$) to claim that $f-f_m=f-diag(A^m)$ is of the form $A^0+diag(A^{1})\cdots+diag(A^{m-1})$ for certain $k$-additive functions $A^k:S^k\to H$, $k=0,1,\cdots, m-1$. 

The uniqueness of the representation \eqref{rep} also follows by induction using the polarization formula.  Indeed, if 
\[
f(x)= A^0+diag(A^{1})(x)\cdots+diag(A^{m})(x)=B^0+diag(B^{1})(x)\cdots+diag(B^{m})(x) \text{ for all } x\in G, 
\]
where $A^{k}, B^{k}:G^k\to H$ are $k$-additive symmetric functions, $k=0,1,\cdots,m$, then 
\[
\Delta_{x_1}\Delta_{x_2}\cdots \Delta_{x_m} f(x)= m!A^m(x_1,\cdots,x_m)=m!B(x_1,\cdots,x_m) \text{ for all } x_1,\cdots,x_m\in G. 
\]
Hence $A^m=B^m$ and we can apply the induction step to the polynomial function  
\[
f(x)-diag(A^{m})(x)= A^0+diag(A^{1})(x)\cdots+diag(A^{m-1})(x)=B^0+diag(B^{1})(x)\cdots+diag(B^{m-1})(x).
\]
Hence $A^k=B^k$ for $k=0,1,\cdots, m$.  This ends the proof. {\hfill $\Box$}

Now we can use a recent result by Kurcharski and Lukasik \cite[Theorem 6]{KuLu} to prove that polynomial functions defined on cancellative semigroups can (uniquely) be extended to polynomial functions defined on the natural extension of their domain. 
\begin{theorem}[Kurcharski and Lukasik] \label{KurLuk}
Let $S$ be a cancellative abelian semigroup, $H$ be an abelian group, $n\in \mathbb{N}$. Furthermore,  let $A_n:S^n\to H$ be an $n$-additive symmetric function. Then, for any abelian group $G$ such that $S\leq G$ and $G=S-S$, the function $A_n$ can be uniquely extended to an $n$-additive symmetric mapping of $G^n$ into $H$. That is, there exists a unique $n$-additive symmetric function $\mathcal{A}_n:G^n\to H$ such that 
\[
\mathcal{A}_n(x_1,\cdots,x_n)=A(x_1,\cdots,x_n) \text{ for all } x_1,\cdots, x_n\in S
\]
\end{theorem}

\begin{corollary}  \label{cor1} Let $S$ be a commutative cancellative semigroup and $H$ be a commutative group. Let $G=S-S$ be a natural extension of $S$. Assume that multiplication by $m!$ is biyective on $H$. Let $f:S\to H$ be a polynomial function of degree at most $m$. Then there exists a unique polynomial function of degree at most $m$, $F:G\to H$, such that $F_{|S}=f$.  
\end{corollary}  

\noindent \textbf{Proof.} It is a direct consequence of Theorems \ref{repre}, \ref{KurLuk}. {\hfill $\Box$}

\begin{corollary}\label{cor2}  Let $S$ be a commutative cancellative semigroup and $H$ be a commutative group. Let $G=S-S$ be a natural extension of $S$. Assume that multiplication by $m!$ is biyective on $H$. Let $f:S\to H$ be a solution of Aichinger equation 
\[
f(x_1+\cdots+x_{m+1})=\sum_{i=1}^{m+1}g_i(x_1,x_2,\cdots, \widehat{x_i},\cdots, x_{m+1}), x_1,\cdots, x_{m+1}\in S  
\]
for certain functions $g_i:S^{m}\to H$. Then there exists a polynomial function of degree at most $m$, $F:G\to H$, such that $F_{|S}=f$.  Moreover, the functions $g_i$ are also polynomial functions of degrees at most $m$.  
\end{corollary}  

\noindent \textbf{Proof.} It is a direct consequence of Theorem \ref{main} and Corollary  \ref{cor1}. {\hfill $\Box$}

When we apply Corollary \ref{cor2} to functions defined on $S=\mathbb{R}_+^p$ with values on $H=\mathbb{K}$ (with, e.g., 
$\mathbb{K}=\mathbb{R} $ or $\mathbb{K}=\mathbb{C}$), we conclude that, if $f: \mathbb{R}_+^p\to \mathbb{K}$  satisfies 
\begin{equation}\label{ARp}
f(x_1+\cdots+x_{m+1})=\sum_{i=1}^{m+1}g_i(x_1,x_2,\cdots, \widehat{x_i},\cdots, x_{m+1}) \text{ for all }x_1,\cdots, x_{m+1}\in \mathbb{R}_+^p  
\end{equation}
then $f$ is the restriction to $\mathbb{R}_+^p$ of a unique polynomial function $F:\mathbb{R}^p\to \mathbb{K}$ of degree at most $m$. In particular, we can use the regularity properties of polynomial functions defined on the ordinary Euclidean space  $\mathbb{R}^p$  (see, e.g., \cite{Lreg}, \cite{laszlo1}) to conclude that, if $f$ is bounded on a set $A\subseteq \mathbb{R}_+^p$ with positive Lebesgue measure and is a solution of Aichinger's equation \eqref{ARp}, there exists a unique ordinary polynomial $F$  of $p$ variables and total degree $\leq m$ (which is of course defined on the whole space $\mathbb{R}^p$) such that $F_{|\mathbb{R}_+^p}=f$, and the functions $g_i$ are also restrictions to $\mathbb{R}_+^{pm}$ of ordinary polynomials of total degree $\leq m$ defined on $\mathbb{R}^{pm}$. 

\section{Degree of the composition and product of polynomial functions}

Theorem \ref{repre}, the canonical representation of polynomial functions, can be used to demonstrate, with an alternative proof (different from the ones given in \cite{AM}, \cite{L}) that compositions and products of polynomial functions are again polynomial functions, and to obtain the natural bounds for their degrees.

Concretely, if we assume that $S$ is a commutative semigroup, $G,H$ are commutative groups, 
$f=f_0+\cdots+f_n$, $f:S\to G$ and $g=g_0+\cdots+g_m$, $g:G\to H$ are polynomial functions with $f_i=diag(A^i)$ and $g_j=diag(B^j)$ for all $i,j$ with $A^i:S^i\to G$ and $B^j:G^j\to H$ symmetric multiadditive functions, then
\begin{eqnarray*}
(g\circ f)(x) &=& \sum_{j=0}^mg_j(\sum_{i=0}^nf_i(x))\\
&=& \sum_{j=0}^mB^j(\sum_{i=0}^nf_i(x),\sum_{i=0}^nf_i(x), \cdots, \sum_{i=0}^nf_i(x))\\
&=& \sum_{j=0}^m\sum_{h_1,h_2,\cdots,h_j=0}^nB^j(f_{h_1}(x),f_{h_2}(x), \cdots, f_{h_j}(x))\\
&=& \sum_{j=0}^m\sum_{h_1,h_2,\cdots,h_j=0}^nB^j(A^{h_1}(x,\cdots,x),A^{h_2}(x,\cdots,x), \cdots, A^{h_j}(x,\cdots,x))\\
&=& \sum_{j=0}^m\sum_{h_1,h_2,\cdots,h_j=0}^n diag(C^{h_1,\cdots,h_j}_j)(x)
\end{eqnarray*} 
where 
$$C^{h_1,\cdots,h_j}_j(x_{1,h_1},\cdots,x_{h_j,h_j})=B^j(A^{h_1}(x_{1,h_1},\cdots,x_{h_1,h_1}), \cdots, A^{h_j}(x_{1,h_j},\cdots,x_{h_j,h_j}))$$
is $(h_1+h_2+\cdots+h_j)$-additive. Obviously, it could be the case that $C^{h_1,\cdots,h_j}_j$ is not symmetric, but this is not a problem since, if $d(x)=diag(A)(x)$ for a certain $k$-additive map $A$, then $d=diag(A^{sim})(x)$, where $A^{sim}$ is the symmetrization of $A$:
\[
A^{sim}(x_1,\cdots,x_k)=\frac{1}{k!}\sum_{\sigma\in S_k}A(x_{\sigma(1)},x_{\sigma(2)},\cdots,x_{\sigma(k)}),
\]
which is $k$-additive and symmetric. It follows that $g\circ f$ is a sum of monomials of degree at most $nm$, since the maximum value of the sums $h_1+\cdots+h_j$ with $0\leq h_i\leq n$ and $0\leq j\leq m$ is $nm$. Hence $g\circ f:S\to H$ is a polynomial function of degree at most $nm$. 

With respect to the product of polynomial functions, things are easier since, if $f_i=diag(A^i)$  and $g_j=diag(B^j)$ are monomials defined on $S$ with values in $\mathbb{C}$ (or any other field extension $\mathbb{K}$ of $\mathbb{Q}$), then $f_i(x)g_j(x)=diag(C^{i,j})(x)$, where 
$$C^{i,j}(x_1,\cdots,x_i,y_1,\cdots,y_j)=A^i(x_1,\cdots,x_i)B^j(y_1,\cdots,y_j)$$
is $(i+j)$-additive. Hence, if $f,g:S\to \mathbb{K}$ are polynomial functions , $f=f_0+\cdots+f_n$, $g=g_0+\cdots+g_m$, then 
$$(f\cdot g)(x)=\sum_{k=0}^{n+m}\sum_{i+j=k}f_i(x)g_j(x) = \sum_{k=0}^{n+m}\sum_{i+j=k} diag((C^{i,j})^{sym})(x)$$
is a polynomial function of degree at most $n+m$. Thus, we have proved the following:

\begin{theorem} Assume that $S$ is a commutative semigroup, $G,H$ are commutative groups, and $\mathbb{K}$ is a field extension of $\mathbb{Q}$. Assume that $f, g:S\to \mathbb{K}$ and $\frak{f}:S\to G$, $\frak{g}:G\to H$ are polynomial functions. Assume also that multiplication by $(\deg \frak{f})!$ is biyective on $G$ and multiplication by $(\deg \frak{g})!$ is biyective on $H$. Then $f\cdot g$ and $\frak{g}\circ \frak{f}$ are polynomial functions. Furthermore, $$\deg (f\cdot g)\leq \deg (f)+\deg(g)$$ and $$\deg (\frak{g}\circ \frak{f})\leq \deg \frak{g} \cdot \deg \frak{f}.$$
\end{theorem}

If, for functions $f:S\to \mathbb{K}$ where $S$ is a commutative semigroup and $\mathbb{K}$ is a field, we use $\Delta_h^{m+1}f=0$ as the definition of polynomial function of degree $\leq m$, then we can demonstrate (by induction on $m$) the formula 
\[
\Delta_h^{m}(f\cdot g)(x)=\sum_{i=0}^{m}\binom{m}{i}\Delta_h^i f(x)\cdot \Delta_{h}^{m-i}g(x+ih)
\]
and, as a corollary, we get that $\deg (f\cdot g)\leq \deg (f)+\deg(g)$. Indeed, if $\deg (f)=n$, $\deg(g)=m$, then
\begin{eqnarray*}
\Delta_h^{n+m+1}(f\cdot g)(x) &=& \sum_{i=0}^{n+m+1}\binom{n+m+1}{i}\Delta_h^i f(x)\cdot \Delta_{h}^{n+m+1-i}g(x+ih)=0
\end{eqnarray*}
since $i\leq n$ implies $m+1\leq n+m+1-i$, which implies that all summands in the second member of the equality above vanish. Morevoer, if we use $\frac{1}{n!}\Delta_{h}^{n}f(x)=f(h)$  as a definition of generalized monomial of degree $n$, the same idea can be used to demonstrate the following result:
\begin{theorem}\label{productmonomials}
Let $S,R$ be commutative semigroups and $\mathbb{K}$ be a field. Assume that $f,g:S\to \mathbb{K}$ are generalized monomials of degrees $n,m$, respectively. Then their product is a generalized monomial of degree $n+m$, and the same holds with products of the form $f(x)g(y)$ with $f:S\to \mathbb{K}$ and $g:R\to\mathbb{K}$ generalized monomials. 
\end{theorem}
\noindent \textbf{Proof.} 
A direct computation shows that  
\begin{eqnarray*}
\Delta_h^{n+m}(f\cdot g)(x) &=& \sum_{i=0}^{n+m}\binom{n+m}{i}\Delta_h^i f(x)\cdot \Delta_{h}^{n+m-i}g(x+ih)\\
&=& \sum_{i=0}^{n-1}\binom{n+m}{i}\Delta_h^i f(x)\cdot \Delta_{h}^{n+m-i}g(x+ih)+ \binom{n+m}{n}\Delta_h^n f(x)\cdot \Delta_{h}^{m}g(x+nh)\\
& \ & \hspace{1cm} +\sum_{i=n+1}^{n+m}\binom{n+m}{i}\Delta_h^i f(x)\cdot \Delta_{h}^{n+m-i}g(x+ih)\\
&=&  \binom{n+m}{n}\Delta_h^n f(x)\cdot \Delta_{h}^{m}g(x+nh)\\
&=&  \frac{(n+m)!}{n!m!} n!f(h)\cdot m!g(h) = (n+m)!(f\cdot g)(h),
\end{eqnarray*}
which proves that $f\cdot g$ is a generalized monomial of degree $n+m$. 
For the second claim, it is enough to take into account that $\phi(x,y)=f(x)g(y)$ can be written as $\phi=\phi_1\cdot \phi_2$ where $\phi_1(x,y)=f(x)$ and $\phi_2(x,y)=g(y)$, and that, under this notation, we have that 
\[
\Delta_{(h,k)}^s\phi_1(x,y)=\Delta_h^sf(x) \text{ and } \Delta_{(h,k)}^s\phi_2(x,y)=\Delta_h^sg(y)
\]
and use Leibniz's formula por $\phi_1\cdot \phi_2$: 
\begin{eqnarray*}
\Delta_{(h,k)}^{n+m}(\phi)(x,y) &=& \sum_{i=0}^{n+m}\binom{n+m}{i}\Delta_{(h,k)}^i\phi_1(x,y)\cdot \Delta_{(h,k)}^{n+m-i}\phi_2((x,y)+i(h,k)) \\
&=& \sum_{i=0}^{n+m}\binom{n+m}{i}\Delta_h^i f(x)\cdot \Delta_{k}^{n+m-i}g(x+ik)\\
&=& \binom{n+m}{n}\Delta_h^n f(x)\cdot \Delta_{k}^{m}g(x+nk)\\
&=&  \frac{(n+m)!}{n!m!} n!f(h)\cdot m!g(k) = (n+m)!f(h)g(k)\\
&=&  (n+m)! \phi(h,k).
\end{eqnarray*}
{\hfill $\Box$}

\end{document}